\documentclass[10pt]{amsart}
\usepackage{latexsym}
\usepackage{amsmath}
\usepackage{amssymb}
\usepackage[all]{xy}

\numberwithin{equation}{section}
\newtheorem{thm}{Theorem }[section]

\newtheorem{prop}[thm]{Proposition}
\newtheorem{lem}[thm]{Lemma}
\newtheorem{claim}[thm]{Claim}
\newtheorem{defi}[thm]{Definition}
\newtheorem{cor}[thm]{Corollary}
\newtheorem{exmp}[thm]{Example}
\newtheorem{rem}[thm]{Remark}
\newtheorem{prob}[thm]{Problem}
\newtheorem{qu}[thm]{Question}
\newtheorem{conj}[thm]{Conjecture}
\newcommand{\C}{{\mathbb C}}
\newcommand{\Q}{{\mathbb Q}}
\newcommand{\Z}{{\mathbb Z}}

\begin{document}
\title{A rational realization problem in Gottlieb group
} 
\author{Toshihiro YAMAGUCHI}
\date{}
\address{Faculty of Education, Kochi University, 2-5-1,Kochi,780-8520, JAPAN}
\email{tyamag@kochi-u.ac.jp}
\maketitle

\begin{abstract}
We define the fibre-restricted  Gottlieb group
with respect to a fibration $\xi  :X\to E\to Y$ in CW complexes.
It is a subgroup of the Gottlieb group of $X$.
When $X$ and $E$ are finite simply connected, 
its  rationalized  model   is   given by the arguments of derivations of 
Sullivan models 
based on F\'{e}lix, Lupton and Smith \cite{FLS}.
We 
consider the realization problem of groups in a Gottlieb group as fibre-restricted Gottlieb groups
in rational homotoy theory.
Especially we define an invariant named as (Gottlieb) depth of $X$ over $Y$.
In particular, when $Y=BS^1$, it is related to the rational toral rank of $X$.
\end{abstract}

\section{Introduction}
For a fibration   $p:E\to Y$ of  CW complexes with fiber $X$, $X\overset{j}{\to} E\overset{p}{\to} Y$,
recall a Gottlieb's  problem: {\it Which homotopy equivalences of $X$ into itself can be extended to fibre homotopy equivalences 
of $E$ into itself ?} (\cite[\S 5]{G2}).
Let  $aut(p)$ denote the space of unpointed fibre-homotopy self-equivalences of $p$,
which  is the subspace of $autE$
 with $g:E\to E$ satisfying $p\circ g=p$
and $aut(X)$  the space of unpointed homotopy self-equivalences of $X$.
In this paper, we consider 
  the restriction map $$R:aut(p)\to aut(X)$$ from  the  viewpoint  of
 Sullivan minimal model theory.

The above Gottlieb's  problem is homotopically interpreted to estimating the image 
(it is denoted as ${\mathcal F}(E)$ in \cite[page.52]{G2}) of the map $\pi_0(R):\pi_0(aut(p))\to \pi_0(aut(X))$.
Note that Gottlieb showed in \cite[Theorem 3]{G2} that $G_1(X)\cong {\mathcal F}(E_{\infty})$ 
for the universal fibration $X\to E_{\infty}\to B_{\infty}$ with fibre $X$.
Also refer  \cite{Y} and \cite{HYY} as a certain case.
Denote  $aut_1(p)$ and $aut_1(X)$ the identity components.
We propose  as a higher dimensional homotopical version, when $X$ and $E$ are simply connected finite CW complexes, 

\begin{prob}\label{a}
Estimate the image of the map $\pi_n(R):\pi_n(aut_1(p))\to \pi_n(aut_1(X))$ in terms of 
$\pi_n(X)$ for $n>0$.
\end{prob}

We note that the path components of $aut(p)$ are all of the same homotopy type \cite[\S 1]{FLS}.
Recall the n-th {\it Gottlieb group} $G_n(X)$ of a CW complex $X$ for $n>0$, which is 
the subgroup of the 
$\pi_n(X)$ 
consisting of homotopy classes of based maps $a:S^n\to X$ such that
the wedge $(a|id_X):S^n\vee X\to X$ extends
to a map $F_a:S^n\times X\to X$
in the homotpy commutative diagram:
$$ \xymatrix{X\times S^n\ar@{.>}[r]^{{F_a}}& X\\
X\vee S^n\ar[u]^{inc.}\ar[r]^{id_X\vee a}&  X\vee X.\ar[u]_{\nabla}
}$$
By  adjointness, it   is equal to $Im (\pi_n(ev):\pi_n(aut_1X)\to \pi_n(X))$ for the ordinal evaluation map $ev:aut_1X\to X$.
\cite{G1}.
We denote the graded Gottlieb group of $X$ as $G_*(X):=\oplus_{n>0}G_n(X)$.
When $X$ is simply connected and finite LS category,
 $G_{even}(X)\otimes \Q=0$  \cite{FH}.
 In the following,
 we often use  the symbol $G_{*}(X)_{\Q}$ as $G_{*}(X)\otimes \Q$. 
 When $X$ is simply connected and finite complex,  $G_{*}(X)_{\Q}=G_*(X_{\Q})$ 
and $\dim G_*(X)_{\Q}\leq cat_0(X)$  \cite{L2}.
Futhermore, Felix and Halperin 
\cite[page.35]{FH} conjecture that
 $G_n(X)_{\Q} =0$ for all $n\geq 2q$
if $X$ be a finite complex of dimension $q$.
Gottlieb \cite{G2} obtained 
$$ G_n(X)=\bigcup_{\xi} {\rm Im} \ \partial^{\xi}_{n+1}$$
for the homotopy connecting homomorphisms $\partial^{\xi}_{n+1}:\pi_{n+1}(Y)\to \pi_n(X)$ of all fibrations $\xi
:X\to E\to Y$. 
A Gottlieb group  contains the following group as a natural subgroup: 


\begin{defi}
Let   $\xi: X\to E \overset{p}{\to} Y$ ba a fibration. 
Then  the n-th fibre-restricted  Gottlieb group of $X$
with respect to $\xi$,
denoted by  
$G^{\xi}_n(X)$  is defined as ${\rm Im} \ \pi_n(ev_X\circ R)$
for the evaluation map $ev_X\circ R:aut_1 (p)\to aut_1X\to X$.
\end{defi}

\begin{thm}\label{Im}
 For any fibration $\xi: X\to E \to Y$, it holds that 
$$ {\rm Im} \ \partial^{\xi}_{n+1}\subset G^{\xi}_n(X)\subset G_n(X)$$
for all $n>0$. 
\end{thm}


For a map $f:X\to Y$,
$G_n(Y,X;f)$ is defined as the image of the evaluation map
$\pi_n(ev):\pi_n(map_f(X,Y))\to \pi_n(Y)$ \cite{WL}.
Here $map_f(X,Y)$ is the commponent of 
$f:X\to Y$ of  the space of maps from $X$ to $Y$.  

\begin{lem}\label{trivial}For any fibration $\xi: X\to E \overset{p}{\to} Y$,

 (1) When $\xi$ is fibre-trivial;i.e.,
$p$ is the projection $E\simeq X\times Y\to Y$, then  $ {\rm Im} \ \partial^{\xi}_{n+1}=0$ and $G^{\xi}_n(X)= G_n(X)$.




(2)  For a fibration $\xi: X\overset{j}{\to} E \to Y$, we have the commutative digram
$$ \xymatrix{G_n^{\xi}(X)\ar[d]_{j_{\sharp}}\ar[r]^{\subset}& G_n(X)\ar[d]^{j_{\sharp}}\\
 G_n(E)\ar[r]^{\subset\ \ \ }& G_n(E,X;j). 
}
$$


\end{lem}


\begin{prop}\label{func}
For $f:Y'\to Y$,
put $f^*\xi :X\to E'\to Y'$ the pull-back fibration of $f$.
Then there is an inclusion $G^{\xi}_n(X)\subset  G^{f^*\xi}_n(X)$ for all $n>0$.
\end{prop}

In this paper, we use the \textit{Sullivan minimal model} $M(Y)$ 
 of a simply connected  space $Y$ of finite type.
It is a free $\Q$-commutative differential graded algebra (DGA) 
 $(\Lambda{V},d)$
 with a $\Q$-graded vector space $V=\bigoplus_{i> 1}V^i$
 where $\dim V^i<\infty$ and a decomposable differential,
 i.e., $d(V^i) \subset (\Lambda^+{V} \cdot \Lambda^+{V})^{i+1}$
 and $d \circ d=0$.
Here  $\Lambda^+{V}$ is 
 the ideal of $\Lambda{V}$ generated by elements of positive degree.
 We often denote  $(\Lambda{V},d)$ simply as $\Lambda V$. 
Denote the degree of an element $x$ of a graded algebra as $|{x}|$. 
Then  $xy=(-1)^{|{x}||{y}|}yx$ and $d(xy)=d(x)y+(-1)^{|{x}|}xd(y)$. 
Notice that $M(Y)$ determines the rational homotopy type $Y_{\Q}$ of $Y$.
In particular, there is an isomorphism $Hom_i (V,\Q )\cong \pi_i(Y)\otimes {\Q}$.
See \cite{FHT} for a general introduction and the standard notations.
 
Let $A$ be a DGA $A=(A^*,d_A)$ with $A^*=\oplus_{i\geq 0}A^i,\ A^0=\Q$, $A^1=0$  and 
the augmentation $\epsilon:A\to \Q$.
Define  
 $Der_i A$  the vector space of derivations of $A$
decreasing the degree by $i>0$,
where  $\theta(xy)=\theta(x)y+(-1)^{i|x|}x\theta(y)$
for $\theta\in Der_iA$. 
We denote  $\oplus_{i>0} Der_iA$ by
$DerA$.
The boundary operator $\delta : Der_* A\to Der_{*-1} A$
is 
defined by $\delta (\sigma)=d_A\circ \sigma-(-1)^{|\sigma |}\sigma\circ d_A$.
It is known  that 
\begin{thm}\cite[\S 11]{Su}(\cite{Sa})\label{Sul}
$\pi_*(aut_1( Y))_{\Q}\cong H_*(Der M(Y))$.
\end{thm}
It  is the case that the base is a point in Theorem \ref{FLS}.
In this paper, according to (\cite[p.314]{Su}),
the symbol $(v,h)\in Der_{|v|-|h|} (\Lambda V)$
 means the derivation sending an element $v\in V$
to $h\in \Lambda V$ and the other to zero.
In particular  $(v,1)=v^*$.
The differential is given as 
$\delta(v,h)\ =\
d\circ(v,h)-(-1)^{|v|-|h|}(v,h)\circ d$.


\begin{thm}\label{FHG} {\rm (\cite{FH})}
For the minimal model $M(X)=(\Lambda W,d)$ of a simply connected finite complex
$X$ and the argumentation $\epsilon:\Lambda W\to \Q$,
$$G_n(X)_{\Q}\cong Im \ (H_n(\epsilon_*):H_n(Der(\Lambda W,d))\to Hom_n(W,\Q)=Hom(W^n,\Q ))$$
for $n>0$. 
\end{thm}



 
Recently,   Y.F\'{e}lix, G.Lupton and S.B.Smith showed  

\begin{thm}\label{FLS}  {\rm  (\cite[Theorem 1.1]{FLS})}
For a   fibration  $p:E\to Y$ with $E$ and $Y$ simply connected and $E$ finite,
 $$\pi_n(aut_1(p))\otimes \Q\cong H_n(Der_{\Lambda V}(\Lambda V\otimes \Lambda W))$$
 for $n>0$.
 \end{thm}
\noindent
Here  $\Lambda V\to \Lambda V\otimes \Lambda W$
is the Koszul-Sullivan model of $p:E\to Y$
 and $(Der_{\Lambda V}(\Lambda V\otimes \Lambda W),\delta)$
is the chain complex of derivations vanishing on  $\Lambda V$.
 Under a general study of  \cite[\S 3]{FLS}
in rational model, we have 

\begin{thm}\label{Y}When $E$ is finite,
the map $\pi_n(R)_{\Q}:\pi_n(aut_1(p))_{ \Q}\to \pi_n(aut_1(X))_{ \Q}$ 
is given by
$$H_n(res):H_n(Der_{\Lambda V}(\Lambda V\otimes \Lambda W),\delta )\to H_n(Der( \Lambda W))$$induced by the restriction of derivations for $n>0$.
Moreover, when $X$ is finite, 
$$G^{\xi}_n(X)_{\Q}=Im \big( H_n(\epsilon_*)\circ H_n(res):H_n(Der_{\Lambda V}(\Lambda V\otimes \Lambda W))\to Hom (W^n,\Q )\big)$$
for the  argumentation $\epsilon :\Lambda W\to \Q$.

\end{thm}

\begin{lem}
When $X$ is an elliptic finite complex, for any fibration  $\xi
:X\to E\to Y$ with $E$ finite,
$$G^{\xi}_N(X)_{\Q} = G_N(X)_{\Q} =\pi_N(X)_{ \Q}$$
where $N:=\max \{ n \mid \pi_n(X)_{ \Q}\neq  0\}$.
Thus 
 $G^{\xi}_*(X)_{ \Q}\neq 0$ when $X$ is elliptic.
\end{lem}








In this paper, we are interested in the following

\begin{prob}
For which  subgroup $G$ of $G_n(X)$
does (not) there  exist  a fibration $\xi$ such that $G^{\xi}_n(X)=G$ ?   
\end{prob}


The following rational invariant is a measure
for a variety of such subgroups. 

\begin{defi}
The rational Gottlieb depth of $X$ with respect to  $Y$ is given as   
$${\rm depth}_Y(X)=\max \bigg\{ n \mid G_*(X)_{\Q}\supset G^{\xi_0}_*(X)_{ \Q}\supsetneq G^{\xi_1}_*(X)_{\Q}\supsetneq\cdots\supsetneq  G^{\xi_n}_*(X)_{ \Q} \bigg\}
$$for some  fibrations $\xi_i : X\to E_i\to Y$ with the  total spaces  $E_i$ finite.
Here $G^{\xi_n}_*(X)_{ \Q}$ may be zero.
If there does  not exist such a fibration $\xi_i$,
let  ${\rm depth}_Y(X)=-1$.
It may occur when $Y$ is infinite.
Also denote the rational Gottlieb depth of $X$ as $${\rm depth}(X)=\max \big\{ {\rm depth}_Y(X)\mid Y \mbox{ is simply connected} \big\}.$$
\end{defi}

From the definition, we have $0\leq {\rm depth}(X)\leq \dim G_*(X)_{ \Q}$.
For a  Lie group $X$, we have  
${\rm depth}(X)< {\rm rank} X$. 
For two fibrations $\xi_i: X_i\to E_i\to Y_i$ ($i=1,2$),
there is the product fibration $\xi_1\times \xi_2:X_1\times X_2\to E_1\times E_2\to Y_1\times Y_2$.
Thus

\begin{lem}For any simply connected space $Y$, it holds that
$${\rm depth}_{Y_1\times Y_2}(X_1\times X_2)\geq {\rm depth}_{Y_1}(X_1)+{\rm depth}_{Y_2}(X_2).$$
Therefore
${\rm depth}(X_1\times X_2)\geq {\rm depth}(X_1)+{\rm depth}(X_2).$
\end{lem}




\begin{prob}
When is the Gottlieb  depth of a space zero ?
\end{prob}

A space $X$ is said to be  an $F_0$-space if
$H^*(X;\Q)\cong \Q [x_1,..,x_n]/(g_1,..,g_n)$ for some $n$
where  $|x_i|$ are even and $g_1,..,g_n$ is a  regular sequence;i.e.,
each $g_i$ is not a zero divisor in the factor-ring $\Q [x_1,..,x_n]/(g_1,..,g_{i-1})$.\\

Recall the Halperin's
\begin{conj}{\rm \cite{Ha}}
 When $X$ is  an $F_0$-space, 
for any fibration $X\overset{j}{\to} E\to Y$, 
$j^*:H^*(E;\Q )\to H^*(X;\Q )$ is surjective.
\end{conj}

It is equivalent to that 
the Serre spectral sequence
rationally collapses in $E_2$-term.
Refer \cite{FHT},\cite{Lu},\cite{M},\cite{Th} for some other equivalent conditions and some  conditions
satisfying the conjectures.

\begin{thm}\label{Hal}Suppose that the Halperin conjecture is ture.
Then 
${\rm depth}(X)=0$
for any $F_0$-space $X$.
\end{thm}

For compact connected Lie groups $G$ and $H$ where 
$H$ is a subgroup of $G$,
when ${\rm rank}G={\rm rank}H$,
the homogeneous space $G/H$ satisfies the Halperin conjecture \cite{ST}.

\begin{cor}
When ${\rm rank}G={\rm rank}H$, it holds that 
${\rm depth}(G/H)=0$.
\end{cor}

Note that the homogeneous space $X=SU(6)/S(3)\times SU(3)$
is not an $F_0$-space.
But ${\rm depth}(X)=0$
by degree arguments. 

A space $X$ is said to be {\it  formal}
if there exists  a quasi-isomorphism $M(X)\to (H^*(X;\Q ),0)$.
A space $X$ is said to be {\it  elliptic}
if its rational homotopy group and rational homology group are both finite.
For  an elliptic space $X$,
$\chi_{\pi}(X):=\dim \pi_{even}(X)_{\Q}-\dim \pi_{odd}(X)_{\Q}$
is said as the {\it homotopy Euler number} of $X$.
When $\dim H^*(X;\Q)<\infty$,
$X$ is an $F_0$-space if and only if $\chi_{\pi}(X)=0$.
\begin{cor}\label{Fel}
When { Halperin} conjecture is ture,
for a formal elliptic space $X$,
${\rm depth}(X)\leq \dim G_*(X)_{\Q}-\dim \pi_{even}(X)_{\Q}$.
\end{cor}

The right hand of this inequation is less than or equal to $ - \chi_{\pi}(X)$.

\begin{qu}
When $X$ is elliptic, does it hold that ${\rm depth}(X)\leq - \chi_{\pi}(X)$ ?
\end{qu}


Let $T^r$ be  an $r$-torus
 $S^1 \times\dots\times S^1$($r$-factors)
and let $r_0(X)$ be the {\it rational  toral rank}, which is 
 the largest integer $r$ such that a
 $T^r$  can act continuously
 on a CW-complex  in  the rational homotopy type of $X$
 with all its isotropy subgroups finite \cite{H}.
Such an action is called  {\it almost free}. 
Refer \cite{JL} for a  relation with  
Gottlieb groups and rational toral ranks.
It is known that $r_0(X)\leq - \chi_{\pi}(X)$ for an elliptic space $X$ \cite{AH}.
The following 
 criterion of Halperin
is  used in this paper.
\begin{prop}\cite[Proposition 4.2]{H}\label{H}
Suppose that $X$ is a simply connected CW-complex  with 
$\dim H^*(X;\Q)<\infty$.
Denote  $M(X)=(\Lambda V,d)$.
Then  $r_0(X) \ge r$ if and only if there is a relative  model $$
(\Q[t_1,\dots,t_r],0)
 \to (\Q[t_1,\dots,t_r] \otimes \Lambda {V},D)
 \to (\Lambda {V},d)$$
where    with $|{t_i}|=2$  and
$Dv \equiv dv$ modulo the ideal $(t_1,\dots,t_r)$ 
 satisfying $\dim H^*(\Q[t_1,\dots,t_r] \otimes \Lambda {V},D)<\infty$.
 Moreover,
 if  $r_0(X) \ge r$,
 then $T^r$ acts freely on a finite complex $X'$
 in the  rational homotopy type of  $X$
 and $M(ET^r\times_{T^r}X')=M(X'/T^r)\cong 
(\Q[t_1,\dots,t_r] \otimes \Lambda {V},D)$.
\end{prop}

\begin{cor}
(1) $r_0(X)=0$ if and only if ${\rm depth}_{BS^1}(X)=-1$.

(2) $r_0(X)=k>0$ if and only if ${\rm depth}_{BT^k}(X)\geq 0$
and ${\rm depth}_{BT^{k+1}}(X)=-1$.
\end{cor}

\begin{thm}\label{Y2}
Suppose that  a torus $T$ acts almost  freely on $X$.
Then
 we have
 $G^{\xi}_n(X)_\Q=  G_n(X/T)_\Q$
for   the Borel  fibration $ X\to ET\times_TX \to BT$.
\end{thm}

For two fibrations $\xi_i:X\to E_i\to Y$ ($i=1,2$),
does there exist a fibration $\xi_3:X\to E_3\to Y$
such that $G^{\xi_1}_*(X)_{\Q}\cap G^{\xi_2}_*(X)_{\Q}=G^{\xi_3}_*(X)_{\Q}$ ?


\begin{defi} Denote  the rational fibre-restrcited  Gottlieb  poset of $X$ with respect to  $Y$ as
the finite poset by inclusions 
$${\mathcal G}_Y(X)=( \{ G^{\xi }_*(X)_{\Q}\}_{\xi_{\Q}},\subset )$$
for   all   fibrations $\xi  : X\to E\to Y$ with the total spaces  $E$ finite. 
\end{defi}



When $Y$ is not finite,  ${\mathcal G}_Y(X)$
may be the empty set.




\begin{thm}\label{hasse}
Suppose that $X\simeq_{\Q} S^{n_1}\times  S^{n_2}\times S^{n_3}\times S^{n_4}$ where $n_i$ is odd
and $Y=\C P^{\infty}$.  
Then the Hasse diagram of ${\mathcal G}_Y(X)$ is one of 

{\small 
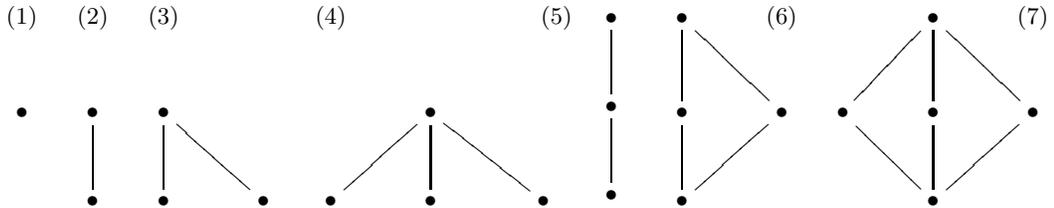
$$
\xymatrix{(1)\\
\bullet
}\ \  \ \ \xymatrix{(2)\\
\bullet\ar@{-}[d]\\
\bullet
} \ \  \ \ \xymatrix{(3)&\\
\bullet\ar@{-}[d]\ar@{-}[dr]&\\
\bullet& \bullet
} \ \ \ \ \ \xymatrix{(4)&&\ \ \ (5)\\
&\bullet\ar@{-}[dl]\ar@{-}[d]\ar@{-}[dr]&\\
\bullet&\bullet& \bullet
}\ \ \ \xymatrix{\bullet\ar@{-}[d]\\
\bullet\ar@{-}[d]\\
\bullet
}\ \ \ \ \ \ \ 
\xymatrix{\bullet\ar@{-}[rd]\ar@{-}[d]&(6)\\
\bullet\ar@{-}[d] &\bullet\ar@{-}[ld]\\
\bullet&
}\ \ \ \ 
\xymatrix{&\bullet\ar@{-}[ld]\ar@{-}[rd]\ar@{-}[d]&(7)\\
\bullet\ar@{-}[rd] &\bullet\ar@{-}[d] &\bullet\ar@{-}[ld]\\
&\bullet&
}
$$
}
and $(1)\sim (7)$ except  $(5)$ are realized as those of ${\mathcal G}_Y(X)$.
\end{thm}

\begin{claim}Then no  3-dimensional
vector space is realized as  
 $G^{\xi }_*(X)_{\Q}$. 
In $(2),(3)$ and $(4)$,  $G^{\xi }_*(X)_{\Q}$ of 
the bottoms are 2-dimensional and in $(5),(6)$ and $(7)$, they are 1-dimensional. 
Since $(5)$ does not exist, if a  1-dimensional $\Q$-space of $G_*(X)_{\Q}$ is realized as $G^{\xi }_*(X)_{\Q}$,
then simultaneously different two or three 2-dimensional  $\Q$-spaces are realized.
\end{claim}

\begin{cor}\label{coro}
When  the circle $S^1$ acts almost freely on a compact connected Lie group $G$ of rank 4,
suppose  that  $\dim G_*(G/S^1)_{\Q}=1$ or
$r_0(G/S^1)=0$.
Then there  exist finite complexes $\{ X_i\}$ ($i=1,2$ or $i=1,2, 3$)
in the rational homotopy type of  $G$ such that

(a) On each  $X_i$,
there is rationally different free
$S^1$-action.
That is, when  $i\neq j$,
there is no map $F$ between rationalized Borel fibrations ${\xi_i}_{\Q}$ and ${\xi_j}_{\Q}$
$$ \xymatrix{ {X_i}_{\Q}\ar[r] \ar@{=}[d]&(ES^1\times_{S^1}X_i)_{\Q}\ar[r]\ar@{.>}[d]^{\simeq}_F&BS^1_{\Q}\ar@{=}[d]\\
 {X_j}_{\Q}\ar[r] &(ES^1\times_{S^1}X_j)_{\Q}\ar[r]&BS^1_{\Q}.
}$$

(b) When $i\neq j$, $G_*^{\xi_i}(X)_{\Q}\neq  G_*^{\xi_j}(X)_{\Q}$ in $G_*(X)_{\Q}$.

(c) For all $i$, $\dim G_*(X_i/S^1)_{\Q}=2$.

(d) For all $i$, $r_0(X_i/S^1)=1$.
\end{cor}

By the arguments of rationa toral poset structure in \cite{Y2},
we have in  Example \ref{toralposet} a  weaker result of Corollary \ref{coro}.  


In this paper, the results of \S 2 follows  due to Gottlieb \cite{G2}
and the results of \S 3 follows  due to F\'{e}lix, Lupton and Smith \cite{FLS}.
Especially, we consider our argument in the case that $Y=BS^1$,
which is  related to almost free circle action on a space in the rational homotopy type of $X$
due to a Halperin's criterion (Proposition \ref{H}).
In this paper, we often use the same symbol $w$ as the dual element $w^*\in Hom(W^n,\Q )\cong \pi_n(X)_{\Q}$ 
for an element $w\in W^n$.

\section{Ordinary properties}

\noindent
{\it Proof of Theorem \ref{Im}.} It is trivial that $G^{\xi}_n(X)\subset G_n(X)$.
We show $ {\rm Im} \ \partial^{\xi}_{n+1}\subset G^{\xi}_n(X)$
by using the arguments in \cite[\S 3-4]{G2}.
Recall that $\xi$ is given by the pull-back diagram of the universal fibration $\xi_{\infty}:E_{\infty}\to B_{\infty}$ with fibre $X$ by a classifying map $h$:
$$ \xymatrix{
X \ar[d]_j\ar@{=}[r]& X\ar[d]^{j_{\infty}}\\
E \ar[d]_p\ar[r]^{\tilde{h}}& E_{\infty}\ar[d]^{p_{\infty}}\\
Y\ar[r]^{h}&B_{\infty},
}$$
where $h$ is is an inclusion and the restriction of $p_{\infty}$ is just $p$ \cite[page.45]{G2}.  

It induces the maps among the following three  fibrations:
$$ \xymatrix{
aut_1(p) \ar[d]\ar[r]^R& aut_1(X)\ar[d]\ar[r]^{ev}&X\ar[d]^{j_{\infty}}\\
aut_1(E) \ar[d]\ar[r]^{\tilde{R}\ \ }& L^*(X,E_{\infty};\tilde{h|_*})\ar[d]\ar[r]^{\ \ \ ev}&E_{\infty}\ar[d]^{p_{\infty}}\\
Y\ar[r]^{h}&B_{\infty}\ar@{=}[r]&B_{\infty}
}$$
where $L^*(X,E_{\infty};\tilde{h|_*})$ means  
(the component of $\tilde{h|_*}$ of) the space of maps from $X$ into $E_{\infty}$
which are homotopy equivalences from $X$  to any fibre of $E_{\infty}$ \cite[page.49]{G2}.  
 Also
 $\tilde{R}$ is defined by $$\tilde{R}(f)=\tilde{h}\circ f\circ j$$ for $f\in aut_1(E).$ 
It is known that  $\pi_i(L^*(X,E_{\infty};\tilde{h|_*}))=0$ for $i>0$. 
Then for $n>0$ 
$$ \xymatrix{\pi_{n+1}(Y)\ar[r]^{\pi_{n+1}(h)}\ar[d]^{\partial}&\pi_{n+1}(B_{\infty})\ar[d]^{\partial}_{\cong}\ar@{=}[r]&
\pi_{n+1}(B_{\infty})\ar[d]^{\partial^{\infty}_{n+1}}\\
\pi_n(aut_1(p))\ar[r]^{\pi_n(R)}&\pi_n(aut_1(X))\ar[r]^{\ \ \pi_n(ev)}&\pi_n(X)
}$$ 
are both commutative.
See \cite[page.50]{G2} for the right hand diagram.
Thus we have $ {\rm Im} \ \partial^{\xi}_{n+1}={\rm Im} ( \partial^{\infty}_{n+1}\circ\pi_{n+1}(h))\subset 
 {\rm Im}(\pi_n(ev)\circ \pi_n(R))=G^{\xi}_n(X)$.
\hfill\qed\\

\noindent
{\it Proof of Lemma \ref{trivial}.}
(1) In this case, since $R:aut_1(p)\to aut_1(X)$ has a section, the map 
 $\pi_n(R):\pi_n(aut_1(p))\to \pi_n(aut_1(X))$ is surjective for all $n$.
 

(2) It follows from the commutative diagram:
$$ \xymatrix{aut_1(p)\ar[d]_{inc.}\ar[r]^{R}& aut_1(X)\ar[d]^{inc.}\ar[r]&X\ar[d]^j\\
aut_1(E)\ar[r]^{res.\ \ }& map_j(X,E)\ar[r]^{\ \ \ \ ev_E} &E,
}
$$ 
where $res.(f):=f\circ j$.\hfill\qed\\

\noindent
{\it Proof of Proposition  \ref{func}.}
The total space of pull-back fibration is given as $E'=\{ (y,x)\in Y\times E \mid f(y)=p(x)\}$.
For $g\in aut_1(p)$, the following diagram 
$${\small  \xymatrix{
E' \ar[dd]_{p'}\ar[rr]^{\tilde{f} }\ar@{.>}[rd]^{f^*(g)}& &E\ar[dd]^{p}\ar[dr]^g&\\
&E' \ar[ld]_{p'}\ar[rr]^{\tilde{f}\ \ \ \ \ \ \ \ \ } &&E\ar[dl]^{p}\\
Y'\ar[rr]^{f}&&Y_{}&
}}$$
is commutative when  ${f}^*(g) (y,x):= (y,g(x))$
Thus there is a map
$$f^*:  aut_1(p)\to aut_1(p')$$
and $$R'\circ f^*=R .$$
Thus we have $\pi_n(ev\circ R'\circ f^*)=\pi_n(ev\circ R) .$
\hfill\qed\\

\section{Sullivan models}

For a  fibration  $X\overset{j}{\to} E\overset{p}{\to} Y$,
recall  that  $[E\times S^n,E]_{\bf o/u}$ is the set of homotopy classes over 
 $1_Y$ and under $1_E$ of maps $F$ 
over 
 $1_Y$ and under $1_E$ \cite[p.379]{FLS};i.e.,
$$ \xymatrix{E\ar@{=}[rd]^{1_E}\ar[d]_{i}&\\
E\times S^n\ar@{.>}[r]^{F}\ar[d]_{{p'}} &E\ar[d]^{p_{}}\\
Y\ar@{=}[r]^{1_Y}& Y_{}
}$$
is a commutative  diagram
in which  $i(x)=(x,*)$ for $x\in E$  and $p'\circ i=p$.
Also a homotopy  over 
 $1_Y$ and under $1_E$ 
is a homotopy map preserving  this commutative diagram \cite[p.379]{FLS}.\\

\noindent
{\it Proof of Theorem \ref{Y}.}  By adjointness,
$\pi_n(R):
\pi_n(aut_1(p))\to \pi_n(aut_1(X))$ is given by the restriction map
$$\Phi :[E\times S^n,E]_{\bf o/u}\to \{ f\in [X\times S^n,X_{}]\mid f\circ i_X\simeq 1_X \}$$
with $ \Phi ([F])= [\overline{F}]:=[F\circ (j\times 1_{S^n}) ]$.  
In concrete terms,  $\Phi$  is represented by   the commutative diagram
$$ \xymatrix{X\times S^n\ar[r]^{\overline{F}}\ar[d]_{j\times 1_{S^n}}&X_{}\ar[d]^{j_{}}\\
E\times S^n\ar[r]^{F}\ar[d]_{{p'}} &E\ar[d]^{p_{}}\\
Y\ar@{=}[r]^{1_Y}& Y_{}
}$$
where  $F\circ i_E=1_E$, $\overline{F}:=F\circ (j\times 1_{S^n})$
and   
$p'\circ i_E=p$.
 The rational model is denoted by the  commutative diagram 
$$ \xymatrix{\Lambda W\otimes \Lambda u/u^2& \Lambda W\ar[l]_{{\overline{F}}^*}\\
\Lambda V\otimes \Lambda W\otimes \Lambda u/u^2\ar[u]^{p_V\otimes 1_u}& \Lambda V\otimes \Lambda W\ar[l]_{\ \ \ \ \ \ F^*}\ar[u]_{p_V}\\
\Lambda V\ar[u]\ar@{=}[r] &\Lambda V,\ar[u]
}$$
where $p_V:\Lambda V\otimes \Lambda W\to  \Lambda W$ is the projection
sending $V$ to zero   and $|u|=n$.
By Theorems \ref{Sul} and \ref{FLS}, $\Phi$
 is given as 
$$\Phi' :H_n(Der_{\Lambda V}(\Lambda V\otimes \Lambda W),\delta)\to H_n(Der( \Lambda W),\delta)$$
induced by the chain map
$$\tilde{\Phi'} :(Der_{\Lambda V}(\Lambda V\otimes \Lambda W),\delta)\to 
(Der( \Lambda W),\delta)$$
with the restriction $\tilde{\Phi'} (\sigma )=\overline{\sigma}:=p_V\circ \sigma$
for $\sigma\in {Der}_{\Lambda V}(\Lambda V\otimes \Lambda W)_n$, 
where 
$F^*=1+(-1)^nu\sigma $
and  ${\overline{F}}^*=1+(-1)^nu\overline{\sigma }$.
\hfill\qed\\

\begin{rem}
Let
$\Lambda V\to \Lambda V\otimes \Lambda W$
be the Koszul-Sullivan model of $p$.
Then there is a DGA-exact sequence
$$0\to (V)\overset{i}{\to}  \Lambda V\otimes \Lambda W \overset{j}{\to} \Lambda W\to 0,$$
where $(V)$ is the DG-ideal  of  the DGA $ \Lambda V\otimes \Lambda W$
generated by $V$.
This induces the exact sequence of chain complexes 
$$0\to Der_{\Lambda V}(\Lambda V\otimes \Lambda W,(V))\overset{i_*}{\to}  
Der_{\Lambda V}(\Lambda V\otimes \Lambda W)\overset{j_*}{\to} Der( \Lambda W)\to 0$$
since $Der( \Lambda W)=Der_{\Lambda V}(\Lambda V\otimes \Lambda W,\Lambda W)$.
Then there is an exact sequence of homologies
$$\cdots \overset{i_n}{\to}  
H_n(Der_{\Lambda V}(\Lambda V\otimes \Lambda W))\overset{j_n}{\to} H_n(Der( \Lambda W))$$
$$\overset{\partial_n}{\to} H_{n-1}(Der_{\Lambda V}(\Lambda V\otimes \Lambda W,(V))\overset{i_{n-1}}{\to}  
H_{n-1}(Der_{\Lambda V}(\Lambda V\otimes \Lambda W))\overset{j_{n-1}}{\to}\cdots ,
$$
in which $j_n$ is equal to the map $\Phi' $ in the above proof.
\end{rem}


\noindent
{\it Proof of Theorem \ref{Hal}}.
Let
$H^*(X;\Q)=\Q [x_1,..,x_n]/(g_1,..,g_n)$ and then 
 $M(X)=(\Lambda (x_1,..,x_n,y_1,..,y_n),d)$
with $|x_i|$ even, $|y_i|$ odd, $dx_i=0$ and $dy_i=g_i\neq 0\in \Q [x_1,..,x_n]$.
Note $G_*(X)_{\Q}=\Q (y_1,..,y_n)$.
We suppose $|y_1|\leq \cdots \leq |y_n|$.

Assume ${\rm depth}(X)>0$.
Then $G_*(X)_{\Q}\supsetneq G^{\xi}_*(X)_{\Q}$
for a fibration $\xi:X\overset{j}{\to} E\to Y$, i.e.,  
 $\Q (y_1,..,y_n)\supsetneq G^{\xi }_*(X)_{\Q}$.
Then, from $D\circ D=0$,  there is an element $y_l$   with 
$D(y_l)\notin  \Q [x_1,..,x_n]$
and 
$$D(y_l)=g_l+\sum_{k<l}y_kf_k+\alpha$$ for some non-zero element 
$f_k\in \Q[x_1,..,x_n]\otimes M(Y)^+$
and $\alpha$ is an element  in the ideal $(y_1y_2,y_1y_3,\cdots ,y_{n-1}y_n)$
in $ \Q [x_1,..,x_n]\otimes \Lambda (y_1,..,y_n) \otimes M(Y)$.
We suppose that  $l$ is the minimal index 
of such elements. 
 Since $g_1,..,g_{l-1}$ is a regular sequence, 
we have $a_1g_1+\cdots +a_{l-1}g_{l-1}\neq 0$ for
any $a_1,..,a_{l-1}\in\Q [x_1,..,x_n]$
with at least one $a_i$ non-zero.
Thus when we denote   $D(y_k)=g_k+h_k$ for $k<l$ with 
$$h_k\in  \Q [x_1,..,x_n]\otimes M(Y)^+\oplus  \Q [x_1,..,x_n]\otimes \Lambda^{>1}(y_1,..,y_n) \otimes M(Y)^+,$$
we see that $\sum_{k<l}(g_k+h_k)f_k$ can  not be zero.
Therefore
$$\beta :=D(\sum_{k<l}y_kf_k)=\sum_{k<l}(g_k+h_k)f_k-\sum_{k<l}y_kD(f_k) $$
is not zero and $\beta+D(\alpha )\neq 0$.
 Thus $D\circ D(y_l)=0$ induces $$D(g_l)=-\beta-D(\alpha )\neq 0,$$
that is, there is an element $x_s$ that $D(x_s)\neq 0$.
Therefore $j^*:H^*(E;\Q )\to H^*(X;\Q )$ is not surjective.
\hfill\qed\\

\noindent
{\it Proof of Corollary \ref{Fel}}.
By \cite{Fe},
$$M(X)=(\Lambda (x_1,..,x_m,y_1,..,y_m,z_1,..,z_n),d)$$
where  $|x_i|$ are even, $|y_i|$ and $|z_i|$ are odd,
$dx_i=dz_i=0$
and  $dy_i=g_i+h_i$.
Here $g_i\in \Q [x_1,..,x_m]$, $g_1,..,g_m$ is a regular sequence and
$h_i\in (z_1,..,z_n)$.     
\hfill\qed\\

\noindent
{\it Proof of Theorem \ref{Y2}.} \
Let 
$T=S^1\times \cdots \times S^1$
(k-times) and  $X\to  ET\times_{T}^{\mu}X\to BT$  the Borel fibration
of  a $T$-action  on $X$.
Then the relative model is given as
 $$
(\Q[t_1,\dots,t_r],0)
 \to (\Q[t_1,\dots,t_r] \otimes \Lambda {W},D)
 \to (\Lambda {W},d).$$
where $M(X)= (\Lambda W,d)$.
There is a  commutative diagram
$$ \xymatrix{
H_n(Der_{\Q[t_1,\dots,t_r]}(\Q[t_1,\dots,t_r]\otimes \Lambda W))\ar[d]_{H_n(i)}\ar[r]^{\ \ \ \ \ \ \ \ \ \ \ \ \ \ \ \ \ \  \Phi'} &H_n(Der( \Lambda W))\ar[r]^{\epsilon_*}
& Hom_n(W,\Q)\ar[d]^{Hom(p_{t},\Q)}\\
H_n(Der(\Q[t_1,\dots,t_r]\otimes \Lambda W))\ar[rr]^{ \epsilon_*}&&Hom_n(\Q(t_1,\dots,t_r)\oplus W,\Q)
}$$
where $i:Der_{\Lambda V}(\Lambda V\otimes \Lambda W)\to Der(\Lambda V\otimes \Lambda W)$ is 
the natural  inclusion of complex
and 
 $p_t$ is the projection sending   $t_1,..,t_r$ to zero.
Since $|t_i|=2$,
$H_n(i)$ is an isomorphism for $n>1$.
The above diagram is commutative and then Theorem \ref{Y} 
it induces the commutaive diagram of homotopy groups 

$$ \xymatrix{
\pi_n(aut_1(p))_{\Q}\ar[d]_{\pi_n(i)_{\Q}}\ar[r]^{\pi_n(R)_{\Q}} &\pi_n(aut_1(X))_{\Q}\ar[r]^{\ \ \pi_n(ev_X)_{\Q}}
& \pi_n(X)_{\Q}\ar[d]^{\pi_n(j)_{\Q}}\\
\pi_n(aut_1(ET\times_TX))_{\Q}\ar[rr]^{\pi_n(ev)_{\Q} }&&\pi_n(ET\times_TX)_{\Q}
}$$
where $i:aut_1(p)\to aut_1(ET\times_TX)$ is the natural  inclusion.
 Since $\pi_n(i)_{\Q}$ is an isomorphism for $n>1$,
 we have the theorem.
\hfill\qed\\

The following lemma follows from the Sullivan model arguments.

\begin{lem}\label{toral}
Let  $\xi': X\overset{j'}{\to} E' \overset{p'}{\to} Y'$
be a fibration with $Y'=\C P^{\infty}$ and $\dim H^*(E';\Q)<\infty$.
 Then for the pull-back   fibration $\xi: X\overset{j}{\to} E \overset{p}{\to} Y$
by   the natural inclusion $i:Y:=\C P^m\to Y'$:
$$ \xymatrix{
X\ar[d]_{j}\ar@{=}[r]& X\ar[d]^{j'}\\
E \ar[d]_{p}\ar[r]& E'\ar[d]^{p'}\\
Y\ar[r]^i&Y',
}$$
we have  $E\simeq_{\Q}E'\times S^{2m+1}$
when  $m$ is sufficiently large.
\end{lem}

\noindent
{\it Proof of Theorem  \ref{hasse}}. 
For the rational type $(n_1,n_2,n_3,n_4)$   of $X$,
we can assume $n_1\leq n_2\leq n_3\leq n_4$ without losing generality.
 Let the set of rational types of $X$
$$I=\big\{ (n_1,n_2,n_3,n_4)\in \Z^{\times 4}\mid  
3\leq n_1\leq n_2\leq n_3\leq n_4\mbox{ are all odd }\big\} $$
 divide into the six cases as follows:
\begin{eqnarray*}
(i)& n_1+n_2>n_4 & \\
(ii)&  n_1+n_2<n_4, & n_1+n_3>n_4 \\
(iii)
  & n_1+n_3<n_4, & n_2+n_3>n_4\\
(iv)&n_1+n_2>n_3,  &  n_2+n_3<n_4\\
(v)& n_1+n_2<n_3, & n_1+n_3<n_4 ,\  \ n_2+n_3>n_4 \\
(vi)& n_1+n_2<n_3,& n_2+n_3<n_4. \\
\end{eqnarray*}

 In general, the  differential $D$ of the relative model of a fibration $\xi :X\to E\to Y=\C P^n$
 is given by $D(w_1)=D(w_2)=0$ and 
\begin{eqnarray*}
D(w_3)&=&aw_1w_2t^m,\\
D(w_4)&=&b_1w_1w_2t^{n_1}+b_2w_1w_3t^{n_2}+b_3w_2w_3t^{n_3}
\end{eqnarray*}
for suitable $m,n_i$  and some $a,b_1,b_2,b_3\in \Q$
which often must  be zero by degree arguments.
When $n$ is sufficiently large, under the conditions $(i)\sim (vi)$, we observe that ${\mathcal G}_Y(X)$ are all different. 
For example, in the case of $(ii)$, ${\mathcal G}_Y(X)=\{ G_*(X)_{\Q} ,G_*^{\xi}(X)_{\Q} \}$ in which $D$ are given as 
\begin{eqnarray*}
G_*(X)_{\Q}:&  & a=b_1=b_2=b_3=0\\
G_*^{\xi}(X)_{\Q}:& &a=b_2=b_3=0,\ b_1\neq 0,\\
&& a\neq 0, \ b_1=b_2=b_3=0\ \mbox{ or }\\
&& a\neq 0,\ b_1\neq 0,\ b_2=b_3=0.
\end{eqnarray*}

We see that $(i),(ii),(iii),(iv),(v)$ and $ (vi)$ correspond
to ${\mathcal G}_Y(X)$ whose  Hasse diagrams are $(1) \ \ \Q (w_1,..,w_4)$,
{\small 
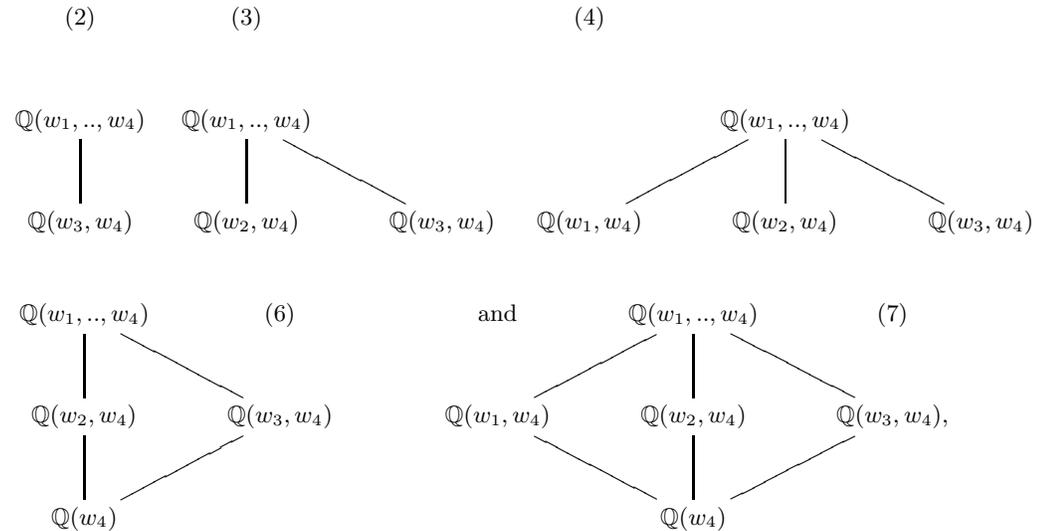
$$
\xymatrix{(2)\\
\Q (w_1,..,w_4)\ar@{-}[d]\\
\Q (w_3,w_4)
} \ \  \ \ \xymatrix{(3)&\\
\Q (w_1,..,w_4)\ar@{-}[d]\ar@{-}[dr]&\\
\Q (w_2,w_4)& \Q (w_3,w_4)
} \ \ \ \ \ \xymatrix{(4)&&  \\
&\Q (w_1,..,w_4)\ar@{-}[dl]\ar@{-}[d]\ar@{-}[dr]&\\
\Q (w_1,w_4)&\Q (w_2,w_4)& \Q (w_3,w_4)
}
$$

$$
\xymatrix{\Q (w_1,..,w_4)\ar@{-}[rd]\ar@{-}[d]&(6)\\
\Q (w_2,w_4)\ar@{-}[d] &\Q (w_3,w_4)\ar@{-}[ld]\\
\Q (w_4)&
}\ \ \ \ \  \ \ \ \ \ \ \ 
\xymatrix{\mbox{and}&\Q (w_1,..,w_4)\ar@{-}[ld]\ar@{-}[rd]\ar@{-}[d]&(7)\\
\Q (w_1,w_4)\ar@{-}[rd] &\Q (w_2,w_4)\ar@{-}[d] &\Q (w_3,w_4),\ar@{-}[ld]\\
&\Q (w_4)&
}
$$
}
\noindent
respectively. For example, they are given by  the table:
\begin{center}
\begin{tabular}{|c||c |c|c|c||c|}
\hline
Hasse diagram&$Dw_1$&$Dw_2$&$Dw_3$&$Dw_4$&basis of $G^{\xi}_*(X)$\\
\hline
$(1)\sim (7)$&$0$&$0$&$0$&$t^*$&$w_1,w_2,w_3,w_4$\\
\hline
$(2)\sim (7)$&$0$&$0$&$0$&$w_1w_2t^*+t^*$&$w_3,w_4$\\
\hline
$(3)\sim (7)$&$0$&$0$&$0$&$w_1w_3t^*+t^*$&$w_2,w_4$\\
\hline
$(4),(7)$&$0$&$0$&$0$&$w_2w_3t^*+t^*$&$w_1,w_4$\\
\hline
$(6),(7)$&$0$&$0$&$w_1w_2t^*$&$w_1w_3t^*+t^*$&$w_4$\\
\hline
\end{tabular}
\end{center}
where $*$ are determined as suitable numbers.
 \hfill\qed\\


\begin{rem}
Note that all above vector spaces of $(1)\sim (7)$ in the above proof are realized as the rational fibre-restricted Gottlieb group of  certain fibrations.
Indeed, for each rational fibration $X\to E\to Y$ with classifying map $h_0$, 
there exists an ordinal classifying map $h$ (with a rational self-equivalence $f$ of $Y$) in the commutative diagram
with rationalization $l_0$:  
$$
\xymatrix{Y\ar@{.>}[rr]^h\ar[d]_{f}&&Baut_1(X)\ar[d]^{l_0}\\
Y\ar[r]^{l_0}&Y_{\Q}\ar[r]^{h_0\ \ \ \ \ }&(Baut_1(X))_{\Q}
}$$
since $Y=\C P^n$ is a simpy connected formal finite complex \cite{Su}(\cite[Remark 3.2]{Pa}).
When $n$ is sufficiently large,
the above arguments applies to the case of  $Y=\C P^{\infty}$
from Lemma \ref{toral}.
\end{rem}






\section{Rational  examples}

\begin{exmp}
For an inclusion $K\to G$ of compact connected simply connected Lie groups,
consider the following fibrations:

(i) $\xi_1:K\to G\to G/K$,

(ii) $\xi_2:G\to G/K\to BK$ and 

(iii) an associated bundle $\xi_3:G/K\to E\to Y$ of a principal $G$-bundle.\\
When $X$ denotes the fibres of them,  we have $G^{\xi_i }_{odd}(X)_{\Q}=\pi_{odd}(X)_{\Q}$ for $i=1,2,3$.
\end{exmp}

\begin{exmp}\label{odd}
If $X$ has the rational homotopy type of the product odd-spheres
$S^{n_1}\times \cdots \times S^{n_k}$  for $k>1$ with $n_1\leq \cdots \leq n_k$,
we have ${\rm depth}(X)=l$ for 
$l=\min  \{ i \mid  n_i=n_{i+1}=\cdots =n_k  \}-1$.
Indeed, for  $m_i=n_k-n_i+1$, let 
$Y=\vee_{i=1}^l S^{m_i}$ be an one point union of $l$ odd-spheres.
Let $M(X)=(\Lambda (w_1,\cdots ,w_k),0)$ with $|w_i|=n_i$
and  $M(Y)=(\Lambda (v_1,\cdots ,v_l,\cdots ),d_Y)$ with with $|v_i|=m_i$
and $d_Yv_1=\cdots =d_Yv_l=0$.
Then we can put 
$M(E)=(\Lambda (v_1,\cdots v_l, w_1,\cdots ,w_k),D)$
where
$Dw_1=\cdots =Dw_N=0$ and 
$$Dw_k=a_1v_1w_1+\cdots +a_lv_lw_l \ ; \ \ \ \ a_i\in \Q$$
When the sequence of $(a_1,..,a_l)$;
$(0,..,0)$, $(1,0,..,0)$, $(1,1,0,..,0), \ \cdots $ and $(1,..,1)$
induces that of inclusions of rational fibre-restricted Gottlibe groups 
respect with $M(E)$
$$\Q (w_1,..,w_k)\supset  \Q (w_2,..,w_k)\supset \Q (w_3,..,w_k)\supset 
\cdots \supset \Q (w_{l+1},..,w_k)$$
In particular, ${\rm depth}(X)=0$ if and only if $n_{1}=\cdots =n_{k}$ for all odd integers.
\end{exmp}

\begin{exmp}\label{hs}
When $X=(S^4\vee S^4)\times S^5$, 
$G_*(S^4\vee S^4)_{\Q}=0$ by \cite{Sm}.
Thus $G_*(X)_{\Q}=G_5(S^5)_{\Q}=\Q$.
Then, by the argument of \cite[Example 6.5]{HS},
there is a fibratuion $\xi :X\to E\to S^2 $ such that 
$G^{\xi}_*(X)_{\Q}=0$.
Thus we have ${\rm depth}_{S^2}(X)=1$.
\end{exmp}

\begin{exmp}
Recall Theorem \ref{Hal}.
Even if $X$ is formal and $\xi:X\overset{j}{\to} E\to Y$ is TNCZ; i.e., $j^*:H^*(E;\Q )\to H^*(X;\Q)$ is surjective, it may holds that $ G^{\xi}_*(X)_{\Q}\neq  G_*(X)_{\Q}$.
Indeed, let $X=S^3\times S^3\times S^4$ where
$M(X)=(\Lambda (w_1,w_2,w_3,w_4),d)$ with
$|w_1|=|w_2|=3$, $|w_3|=4$, $|w_4|=7$, $dw_1=dw_2=dw_3=0$ and
$dw_4=w_3^2$.
Then $G_*(X)_{\Q}=\Q (w_1,w_2, w_4)$.
Let $Y=S^2$ where
$M(Y)=(\Lambda (v_1,v_2),d)$ with
$|v_1|=2$, $|v_2|=3$, $dv_1=0$ and 
$dv_2=v_1^2$.
If the rational model of a fibration $\xi:X\to E\to Y$ is given  as $$Dw_1=Dw_2=Dw_3=0,\ Dw_4=w_1w_2v_1+w_3^2,$$
then on the cohomology of the fibre-inclusion,
$$ j^*:\Lambda (w_1,w_2)\otimes \Q [w_3,v_1]/(w_1w_2v_1+w_3^2,v_1^2)\to 
\Lambda (w_1,w_2)\otimes \Q [w_3]/(w_3^2) $$
is surjective but $G^{\xi}_*(X)_{\Q}=\Q (w_4)$.
Thus ${\rm depth}_{S^2}(X)\geq 1$.


\end{exmp}




\begin{exmp}\label{1}


For  the principal bundle $SU(3)\to SU(5)\to SU(5)/SU(3)$ 
with $SU(n)$ the special unitary group,
the model of the induced fibration 
 $SU(5)\to SU(5)/SU(3)\to BSU(3)$ is given as 
$$(\Q[t_1,t_2],0)\to (\Q[t_1,t_2]\otimes \Lambda (v_1,v_2,v_3,v_4),D)\to (\Lambda (v_1,v_2,v_3,v_4),0)$$
with
$|v_i|=2+2i-1$, $|t_1|=4$, $|t_2|=6$,
$Dv_1=t_1$, $Dv_2=t_2$ and $Dv_3=Dv_4=0$.
Then we have 
 \begin{center}
\begin{tabular}{ |c|c|c|c|c|}
\hline
n&  $\pi_n(aut_1(SU(5)))_{\Q}$&   $\pi_n(aut_1(p))_{\Q}$&$Im\ \pi_n(R)_{\Q}$ & $G^{\xi}_*(SU(5))_{\Q}$\\
\hline
1&$(v_4,v_1v_2)$ &&&\\
\hline
2 &$(v_4,v_3),(v_3,v_2),(v_2,v_1)$  &$(v_4,v_3)$ & $(v_4,v_3)$&\\
\hline
3& $(v_1,1)$& $(v_1,1)$& $(v_1,1)$&$v_1^*$\\
\hline
4& $(v_3,v_1),(v_4,v_2)$&&&\\
\hline
5& $(v_2,1)$& $(v_2,1)$& $(v_2,1)$& $v_2^*$\\
\hline
6& $(v_4,v_1)$&&&\\
\hline
7& $(v_3,1)$& $(v_3,1)$&$(v_3,1)$  &$v_3^*$ \\
\hline
8&  && &\\
\hline
9& $(v_4,1)$&  $(v_4,1)$&$(v_4,1)$  &$v_4^*$\\
\hline
\end{tabular}
\end{center}

Note  $\pi_*(SU(5)/SU(3))_{\Q}=\Q (v_3,v_4)$ and 
 $\pi_*(SU(3))_{\Q}=\Q (v_1,v_2)$.
From the table, 
$G^{\xi}_i(SU(5))_{\Q}=Im\ \pi_i(ev \circ R)_{\Q}=\Q$ for $i=3,5$
 but  $\pi_i(aut_1SU(5)/SU(3))_{\Q}=0$ for $i=3,5$.
 It deduces an example that $G^{\xi}_*(X)_{\Q}=G_*(X/G)_{\Q}$ does not fold in a general free $G$-action on $X$.
Compare with Theorem \ref{Y2}.\\

 




\end{exmp}

\begin{exmp}\label{odd2}
 Let  $X=S^{3  }\times S^{5  }\times S^{9 }\times S^{11  }\times S^{17  }$.  Then 
$\dim G_*(X)_{\Q} =5$.
Denote $M(X)=(\Lambda (w_1,w_2,w_3,w_4,w_5),0)$
of $|w_1|=3$, $|w_2|=5$, $|w_3|=9$, $|w_4|=11$, $|w_5|=17$.
When  $Y=\C P^n$ for $n\geq 5$,  we have the table $(a)$:
\begin{center}
\begin{tabular}{|c|l |}
\hline
dimension &   basis set of ${\mathcal G}_Y(X)$ \\
\hline
5& $\{ w_1,w_2,w_3,w_4,w_5\}$ \\
\hline
4& $\phi$ \\
\hline
3&   $\{ w_1,w_3,w_5\}$, $\{ w_1,w_4,w_5\}$, $\{ w_2,w_3,w_5\}$, $\{ w_2,w_4,w_5\}$, $\{ w_3,w_4,w_5\}$\\
\hline
2&$\{ w_3,w_5\}$, $\{ w_4,w_5\}$ \\
\hline
1&$\{ w_5\}$\\
\hline
0&$\phi$\\
\hline
\end{tabular}
\end{center}
On the other hand, when $Y'=S^2$, we have the table $(b)$:
\begin{center}
\begin{tabular}{|c|l |}
\hline
dimension &   basis set of  ${\mathcal G}_{Y'}(X)$ \\
\hline
5& $\{ w_1,w_2,w_3,w_4,w_5\}$ \\
\hline
4& $\phi$ \\
\hline
3&   $\{ w_1,w_3,w_5\}$, $\{ w_3,w_4,w_5\}$\\
\hline
2&$\{ w_3,w_5\}$ \\
\hline
1&$\phi$\\
\hline
0&$\phi$\\
\hline
\end{tabular}
\end{center}
Thus the Hasse diagrams are given as

{\small 
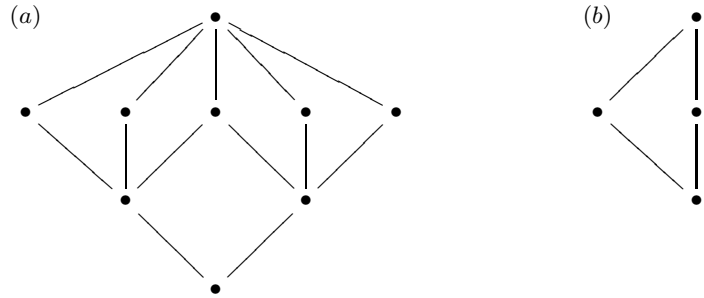
$$
\xymatrix{(a)&&\bullet\ar@{-}[ld]\ar@{-}[lld]\ar@{-}[rd]\ar@{-}[rrd]\ar@{-}[d]&&\\
\bullet\ar@{-}[rd]&\bullet\ar@{-}[d] &\bullet\ar@{-}[rd] \ar@{-}[ld]&\bullet\ar@{-}[d]&\bullet\ar@{-}[ld]\\
&\bullet\ar@{-}[rd]&&\bullet\ar@{-}[ld]&\\
&&\bullet&&
}\ \ \ \ \ \ \ \ \ \ \ \ \ \ \ \ \ \ \ 
\xymatrix{(b) & \bullet\ar@{-}[ld]\ar@{-}[d]\\
\bullet \ar@{-}[rd]& \bullet\ar@{-}[d]\\
& \bullet\\
& 
}$$
}
Note that ${\rm depth}_Y(X)=3$ and ${\rm depth}_{Y'}(X)=2$.

\begin{exmp}\label{(5)}
When is ${\mathcal G}_Y(X)$ well-ordered ?
 Suppose $X=S^{3}\times S^9 \times \C P^{5}\times S^{17}$ and
$Y=\C P^{\infty}$.
Then $M(X)=(\Lambda (u,w_1,w_2,w_3,w_4),d)$ with 
$|u|=2$, $|w_1|=3$, $|w_2|=9$, $|w_3|=11$, $|w_4|=17$, $dw_1=dw_2=dw_4=0$, $dw_3=u^6$.
For $M(\C P^{\infty})=(\Q [t],0)$ ($|t|=2$),
when the model of $\xi$ is given by $Dw_1=Dw_2=0$,　$Dw_3=u^6$,
$Dw_4=w_1w_2t^3+t^9$, 
we have $G^{\xi}(X)_{\Q}=\Q (w_2,w_4)$.
On the other hand,  when  the model of $\xi$ is given by
$$Dw_1=0,\ \ \ Dw_2=u^4t,\ \ \ Dw_3=u^6,\ \ \ Dw_4=w_1w_3t^2-w_1w_2u^2t+t^9,$$
we have $G^{\xi}(X)_{\Q}=\Q (w_4)$.
By degree reason, $\Q (w_1,w_4)$
and $\Q (w_2,w_4)$
are not realized.
Thus ${\mathcal G}_Y(X)$
is given by the well-ordered set
$$\Q (w_1,w_2,w_3,w_4)\supset \Q (w_3,w_4) \supset \Q (w_4),$$
whose Hasse diagram  is $(5)$ in Theorem  \ref{hasse}.
Note that when $X'=S^{3}\times  S^{9}\times S^{11}\times S^{17}$ and 
$Y=\C P^{\infty}$, it is in the case of $(iii)$ in the proof of Theorem \ref{hasse}.
Thus this  example gives that 
there is a map
$f:X'\to X$ that
$f_{\sharp}:G_*(X')_{\Q}\cong G_*(X)_{\Q}$ but
${\rm depth}_Y(X')=1<2={\rm depth}_Y(X)$.
Futhermore, this example indicates that  
 it may be 
${\rm depth}_{BS^1}(X')< {\rm depth}_{BS^1}(X'/S^1)$
for an almost free $S^1$-action on a space $X'$.
\end{exmp}

\end{exmp}

\begin{exmp}\label{toralposet}
Recall the rational toral rank poset of $X$, denoted as 
${\mathcal T}_0(X)$ \cite{Y2}.
When $X\simeq_{\Q} S^{n_1}\times  S^{n_2}\times S^{n_3}\times S^{n_4}$ where $n_i$ is odd,
the Hasse diagram of ${\mathcal T}_0(X)$ is one of $(a),(b)$ and $(c)$:
{\small $$(a) 
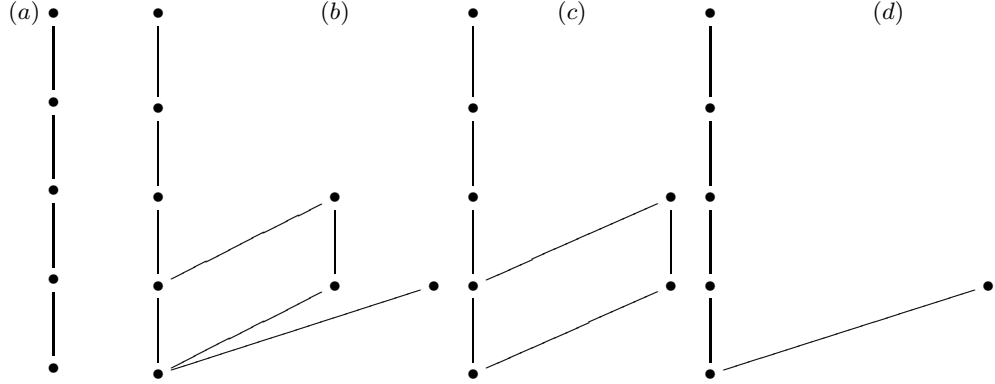
\xymatrix{ 
\bullet &  \\
\bullet \ar@{-}[u]&  \\
\bullet \ar@{-}[u]&   \\
\bullet \ar@{-}[u] &    \\
\bullet\ar@{-}[u]&\\
}\xymatrix{ 
\bullet & & (b) & \\
\bullet \ar@{-}[u]&  &  & \\
\bullet \ar@{-}[u]& & \bullet& \\
\bullet \ar@{-}[u]\ar@{-}[urr] &  & \bullet\ar@{-}[u]& \bullet  \\
\bullet\ar@{-}[u]\ar@{-}[urr]\ar@{-}[urrr]\\
} \ \ \xymatrix{ 
\bullet &(c) &  \\
\bullet \ar@{-}[u]&  &   \\
\bullet \ar@{-}[u]& & \bullet \\
\bullet \ar@{-}[u]\ar@{-}[urr] &  & \bullet\ar@{-}[u]  \\
\bullet\ar@{-}[u]\ar@{-}[urr]\\
}\ \ \xymatrix{ 
\bullet & & (d) & \\
\bullet \ar@{-}[u]&  &  & \\
\bullet \ar@{-}[u]& & & \\
\bullet \ar@{-}[u] &  & & \bullet  \\
\bullet\ar@{-}[u]\ar@{-}[urrr]\\
}$$
}
and $(d)$ is not realized as that of such a space (it corresponds to that
the bottom of $(2)$ in the proof Theorem \ref{hasse} is not given by a fibration $\xi$ with $G_*^{\xi}(X)_{\Q}=\Q (w_4)$).
For $I$ in the proof Theorem \ref{hasse},   
$(a)$ is given when $n_1+n_2>n_4$, 
$(b)$ is given when $n_1+n_2<n_3$, $n_1+n_3<n_4$ and 
$(c)$ is given for the other cases. 
Since $(d)$ is not realized by ${\mathcal T}_0(X) $,  we have the result:
When  the circle $S^1$ acts almost freely on a compact connected Lie group $G$ of rank 4, 
suppose  that $r_0(G/S^1)=0$.
Then there  exist a finite complex $ X$ 
in the same rational homotopy type with $G$ such that
 $r_0(X/S^1)=1$.
(Compare with Corollary \ref{coro}.)
\end{exmp}

\begin{exmp}
When 
$X\simeq_{\Q} S^{n_1}\times  S^{n_2}\times \cdots \times S^{n_m}$ ($n_1\leq n_2\leq \cdots \leq n_m$;odd),
${\rm depth}_{BS^1}(X)=0$
if and only if $n_1+n_2> n_m$. 
Then  $r_0(X/T^r)=r_0(X)-r$ for any $T^r$-action on $X$.
For example, since 
$SU(4)\simeq_{\Q}S^{3}\times  S^{5}\times S^{7}$,
${\rm depth}_{BS^1}(SU(4))=0$.
Then for any almost free $T^r$-action on $SU(4)$, we have
$r_0(SU(4)/T^r)=r_0(SU(4))-r=3-r$.
But for $SU(5)(\simeq_{\Q}S^{3}\times  S^{5}\times S^{7}\times S^{9})$,
${\rm depth}_{BS^1}(SU(5))=1$.
Then there is a  free $S^1$-action on a space $X$ in the rational homotopy type of $SU(5)$ that  $r_0(X/S^1)=1=r_0(SU(5))-3$.
\end{exmp}

\begin{defi} \label{rdepth}
Let the rational toral depth of a simply connected space $X$ 
$$r_0\mbox{-{\rm depth}}(X):= \max \big\{ n \mid (0,1), (r_1,1),   (r_2,1),  \cdots ,(r_n,1) \in {\mathcal T}_0(X) $$
$$\mbox{ \ with\  \  }
0=r_0<r_1<r_2<\cdots <r_n<r_0(X)\big\}.$$
Here we define $r_0\mbox{-{\rm depth}}(X)=-1$ when $r_0(X)=0$.
\end{defi}

Then $r_0\mbox{-depth}(X)=0$ if and only if there is an almost free $S^1$-action on $X'$ and 
 $r_0(X'/S^1)=r_0(X')-1$ for any almost free $S^1$-action on $X'$ in the rational homotopy type of $X$.
For example,  in Example \ref{toralposet},  $r_0\mbox{-depth}(X)$ of $(a),(b)$ and $(c)$ are $0,2$ and $1$, respectively. 
In general, $r_0\mbox{-{\rm depth}}(X)\geq  r_0\mbox{-{\rm depth}}(X/S^1)$
for an almost free $S^1$-action on a space $X$.(Compare with Example \ref{(5)}.)

 \begin{lem}
$r_0\mbox{-{\rm depth}}(X\times Y)\geq r_0\mbox{-{\rm depth}}(X)+r_0\mbox{\rm -depth}(Y)$
 \end{lem}

There is an example that $r_0\mbox{-depth}(X\times S^{12})\geq 0>-2=-1+-1 =r_0\mbox{-depth}(X)+r_0\mbox{-depth}(S^{12})$ due to Halperin
\cite{JL}.
 Recall in Example \ref{toralposet} that $r_0\mbox{\rm -depth}(X)={\rm depth}_{BS^1} (X)$ when $X\simeq_{\Q} S^{n_1}\times  S^{n_2}\times S^{n_3}\times S^{n_4}$
 ($n_i$  odd).
 They are $0,1$ or $2$.
 

\begin{qu} When $X$ is a Lie group, 
$r_0\mbox{-{\rm depth}}(X)=\mbox{\rm depth}_{BS^1}(X)$ ?
\end{qu}

\end{document}